\documentclass[11pt,twoside,letterpaper]{article} 
\usepackage{times,fancyhdr}   
\usepackage{amsfonts}                 
\usepackage{amsthm}                        
\usepackage{amsmath}                         
\usepackage{amssymb}  
  
\numberwithin{equation}{section}


\makeatletter 
\def\cleardoublepage{\clearpage\if@twoside \ifodd\c@page\else%
    \hbox{}%
    \thispagestyle{empty}%
    \newpage%
    \if@twocolumn\hbox{}\newpage\fi\fi\fi} 
\makeatother 

\def \N{\mathbb{N}}

\def \Q{\mathbb{Q}}

\def \Z{\mathbb{Z}}

\def \a{\alpha}

\def \E{\mathcal{E}}
\def \EE{\mathcal{E}^*}

\newtheorem{thm}{Theorem}[section] 

\newtheorem{cor}{Corollary}[section]
\newtheorem{lemma}{Lemma}[section]

\begin{document}
\title{
{\begin{flushleft}
\vskip 0.45in
{\normalsize\bfseries\textit{ }}
\end{flushleft}
\vskip 0.45in
\bfseries\scshape On error sums formed by rational approximations with split denominators}}

\thispagestyle{fancy}
\fancyhead{}
\fancyhead[L]{In: Book Title \\ 
Editor: Editor Name, pp. {\thepage-\pageref{lastpage-01}}} 
\fancyhead[R]{ISBN 0000000000  \\
\copyright~2007 Nova Science Publishers, Inc.}
\fancyfoot{}
\renewcommand{\headrulewidth}{0pt}

\author{\bfseries\itshape Thomas Baruchel and Carsten Elsner\thanks{Lyc{\'e}e naval, Centre d'instruction naval, 
Brest, France, e-mail: baruchel@riseup.net
\newline Fachhochschule f{\"u}r die Wirtschaft, University of Applied Sciences, Freundallee 15,
D-30173 Hannover, Germany, e-mail: carsten.elsner@fhdw.de }}

\date{}
\maketitle
\thispagestyle{empty}
\setcounter{page}{1}

\[\]
\[\]
\begin{abstract}
In this paper we consider error sums of the form
\[\sum_{m=0}^{\infty} \varepsilon_m\Big( \,b_m\alpha - \frac{a_m}{c_m}\,\Big) \,,\]
where $\alpha$ is a real number, $a_m$, $b_m$, $c_m$ are integers, and $\varepsilon_m=1$ or $\varepsilon_m ={(-1)}^m$.
In particular, we investigate such sums for 
\[\alpha \in \big\{ \pi, e,e^{1/2},e^{1/3},\dots, \log (1+t), \zeta(2), \zeta(3) \big\} \]
and exhibit some connections between rational coefficients occurring in error sums for Ap{\'e}ry's continued fraction for
$\zeta(2)$ and well-known integer sequences. The concept of the paper generalizes the theory of ordinary error sums,
which are given by $b_m=q_m$ and $a_m/c_m=p_m$ with the convergents $p_m/q_m$ from the continued fraction expansion of
$\alpha$. 
\end{abstract}
 
\vspace{.08in} \noindent \textbf{Keywords:} Errorr sums, continued fractions   \\
\noindent \textbf{AMS Subject Classification:} 11J70, 11J04, 33B10, 33C05, 11M06. \\


\pagestyle{fancy}  
\fancyhead{}
\fancyhead[EC]{Thomas Baruchel and Carsten Elsner}
\fancyhead[EL,OR]{\thepage}
\fancyhead[OC]{On error sums formed by rational approximations with split denominators}
\fancyfoot{}
\renewcommand\headrulewidth{0.5pt} 

\newpage
\section{Introduction} \label{Sec1}   
Let $\a$ be a real number. We assume that there is a  sequence $B:={(b_n)}_{n\geq 0}$ of integers, a sequence
$R:={(r_n)}_{n\geq 0}$ of rationals $r_n=a_n/c_n$, say, with $a_n\in {\Z}$ and $c_n\in {\N}$, and a real number $\omega >1$
satisfying
\begin{equation}
\big| b_nc_n\a - a_n \big| \,\ll \, \frac{c_n}{\omega^n} \qquad (n\geq 0) \,.
\label{10}
\end{equation}
This is equivalent with
\begin{equation}
\big| b_n\a - r_n \big| \,=\,  \Big| b_n\a - \frac{a_n}{c_n} \Big|  \, \ll \, \frac{1}{\omega^n} \qquad (n\geq 0) \,.
\label{20}
\end{equation}
We consider the fraction $a_n/b_nc_n$ as a rational approximation of $\a$ with split denominator $b_nc_n$. 
Since $\omega >1$, the error sums
\begin{eqnarray}
\EE \big( B,R,\a \big) &:=& \sum_{m=0}^{\infty} \big( b_m\a - r_m \big) \,=\, \sum_{m=0}^{\infty} \Big( b_m\a - 
\frac{a_m}{c_m} \Big) \,,
\label{30} \\
\E \big( B,R,\a \big) &:=& \sum_{m=0}^{\infty} \big| b_m\a - r_m \big| \,=\, \sum_{m=0}^{\infty} \Big| b_m\a - 
\frac{a_m}{c_m} \Big| 
\label{40} 
\end{eqnarray}
exist. Let ${(p_n/q_n)}_{n\geq 0}$ be the sequence of convergents of $\a$ defined by $p_n/q_n = \langle a_0;a_1,a_2,\dots a_n
\rangle$ from the regular continued fraction expansion
\[\a \,=\, \langle a_0;a_1,a_2,\dots \rangle \,=\, a_0 + \cfrac{1}{a_1+\cfrac{1}{a_2+\cfrac{1}{\ddots}}} \]
of $\a$. The error sums of $\a$ for $B={(q_n)}_{n\geq 0}$ and $R={(p_n)}_{n\geq 0}$, namely
\begin{eqnarray*}
\EE(\a) &:=& \EE \big( B,R,\a \big) \,=\, \sum_{m=0}^{\infty} \big( q_m\a - p_m \big)  \,,\\
\E(\a)  &:=& \E \big( B,R,\a \big) \,=\, \sum_{m=0}^{\infty} \big| q_m\a - p_m \big| \,,
\end{eqnarray*}   
were already studied in some papers \cite{ElsnerA, ElsnerB, ElsnerC, ElsnerD}. We call $\EE(\a)$ and $\E(\a)$ {\em ordinary error sums\/}. 
Conversely, for $B={(1)}_{n\geq 0}$ and $R={(p_n/q_n)}_{n\geq 0}$, until now nobody has found any remarkable approach to the error sums
\begin{eqnarray*}
\EE(B,R,\a) &=& \sum_{m=0}^{\infty} \Big( \a - \frac{p_m}{q_m} \Big)  \,,\\
\E(B,R,\a)  &=&  \sum_{m=0}^{\infty} \Big| \a - \frac{p_m}{q_m} \Big|  \,.
\end{eqnarray*}   
In this paper we focus our interest on the series in (\ref{30}) and (\ref{40}) in the case of particular values of $\a$ and well-known
rational approximations of the form
\begin{equation}
0 \,<\, \Big| b_n\a - \frac{a_n}{c_n} \Big| \,\ll \,  \frac{1}{\omega^n} \qquad (n\geq 0) \,.
\label{50}
\end{equation}
Among others we are going to study the numbers
\[\a \,\in \, \Big\{ \,\pi,\, e^{1/l},\,\frac{\log \rho}{\sqrt{5}},\, \log(1+t),\,\zeta(2),\,\zeta(3)\,\Big\} \,,\]
where $l=1,2,\dots$, $e=\exp (1)$, $\rho =(1+\sqrt{5})/2$, and $-1<t\leq 1$, and we shall investigate extraordinary properties of  corresponding error sums (\ref{30}) and
(\ref{40}).

\section{Ordinary error sums for values of the exponential function} \label{Sec2}
Ordinary error sums connected with the exponential function are studied in \cite{Allouche, Elsner1}. Here, our goal is to express this usual error sums itselves by a non-regular
continued fraction. For this purpose we express the error integral
\[\mbox{erf\,}(z) \,=\, \frac{2}{\sqrt{\pi}} \int_0^z e^{-t^2}\,dt \]
by a hypergeometric series, which again can be transformed into a Gauss-type continued fraction. 
\begin{thm}
Let $l\geq 2$ be an integer, and let $p_n/q_n$ denote the convergents of $e^{1/l}$. Then we have
\begin{eqnarray*}
\mathcal{E}(e^{1/l}\,) &=& \sum_{n\geq 0} \big| e^{1/l}q_n - p_n \big| \,=\, e^{1/l} \sqrt{\frac{\pi}{l}} \mbox{erf\,}\Big( \,\frac{1}{\sqrt{l}} \,\Big) \,=\, \frac{2e^{1/l}}{\sqrt{l}} 
\int_0^{1/\sqrt{l}} e^{-t^2} \,dt \\ 
&=& \cfrac{1/l}{1/2 - \cfrac{1/2l}{3/2 + \cfrac{2/2l}{5/2 - \cfrac{3/2l}{7/2+\cfrac{4/2l}{9/2-\cfrac{5/2l}{11/2+\ddots \cfrac{{(-1)}^m m/2l}{(2m+1)/2+\ddots}}}}}}} 
\qquad (m\geq 1) \,.
\end{eqnarray*}
\end{thm}
{\em Proof:\/} \,The first identity of the theorem expressing $\mathcal{E}(e^{1/l})$ by an error integral is already known from
\cite{Allouche, Elsner1}. In order to prove the continued fraction expansion, we set
\[f(z) \,:=\, \frac{\sqrt{\pi}}{2} ze^{z^2} \mbox{erf\,}(z) \,=\, ze^{z^2} \int_0^z e^{-t^2}\,dt \,.\]
We express $f(z)$ in terms of a hypergeometric function ${}_1F_1(\alpha, \beta;z^2)$. 
\begin{eqnarray*}
f(z) &=& ze^{z^2} \int_0^z \sum_{\nu =0}^{\infty} \frac{{(-1)}^{\nu}t^{2\nu}}{\nu!}\,dt \,=\, 
ze^{z^2} \sum_{\nu =0}^{\infty}  \frac{{(-1)}^{\nu}z^{2\nu +1}}{(2\nu +1)\nu!} \\
&=& z^2 \left( \sum_{\mu =0}^{\infty} \frac{z^{2\mu}}{\mu!} \right) \left( \sum_{\nu =0}^{\infty}  
\frac{{(-1)}^{\nu}z^{2\nu}}{(2\nu +1)\nu!} \right) \\
&=& z^2 \sum_{\mu =0}^{\infty} \,\sum_{\nu =0}^{\infty} \frac{{(-1)}^{\nu} z^{2(\nu +\mu)}}{(2\nu +1)\nu! \mu!} 
\,=\, z^2\sum_{k=0}^{\infty} \Big( \mathop{\sum_{\mu =0}^{\infty} \,\sum_{\nu =0}^{\infty}}_{\mu + \nu =k} 
\frac{{(-1)}^{\nu}}{(2\nu +1)\nu! \mu!} \Big) z^{2k} \\  
&=& z^2 \sum_{k=0}^{\infty} \Big(\,\sum_{\nu =0}^k \frac{{(-1)}^{\nu}}{(2\nu +1)\nu! (k-\nu)!} \Big) z^{2k} 
\,=\, z^2 \sum_{k=0}^{\infty} \frac{1}{k!}\Big(\,\sum_{\nu =0}^k \frac{{(-1)}^{\nu}{k \choose \nu}}{2\nu +1} \Big) z^{2k} \,.
\end{eqnarray*}
From \cite[p.\ 68]{Larsen}, Remark\,8.5, we have the following formula (with $k$ replaced by $\nu$ and $n$ replaced by $k$)
\[\frac{1}{d^k {[k]}_k} \sum_{\nu =0}^k \frac{{(-1)}^{\nu}{k \choose \nu}}{c+\nu d} \,=\, \frac{1}{c(c+d)(c+2d)\cdots (c+kd)} \,,\]
where ${[k]}_k =k!$. Setting $c=1$ and $d=2$, it follows that
\[\frac{1}{k!}\,\sum_{\nu =0}^k \frac{{(-1)}^{\nu}{k \choose \nu}}{2\nu +1} \,=\, \frac{2^k}{1\cdot 3\cdot 5 \cdots (2k+1)} 
\,=\, \frac{1}{{(3/2)}_k} \,.\]  
This gives
\[f(z) \,=\, z^2\sum_{k=0}^{\infty} \frac{z^{2k}}{{(3/2)}_k} \,=\, z^2 \sum_{k=0}^{\infty} \frac{{(1)}_k}{k!{(3/2)}_k}  z^{2k}
\,=\, z^2\, {}_1F_1\big( 1,3/2;z^2 \big) \,.\]  
The function ${}_1F_1\big( 1,3/2;z^2 \big)$ can be expressed by a Gauss-type continued fraction. Using formula (8) on page\,123
in \cite{Perron} with $\gamma =3/2$ and $x=z^2$, we have
\[{}_1F_1\big( 1,3/2;z^2 \big) \,=\, \cfrac{1/2}{1/2 - \cfrac{z^2/2}{3/2 + \cfrac{2z^2/2}{5/2 - \cfrac{3z^2/2}{7/2+
\cfrac{4z^2/2}{9/2-\cfrac{5z^2/2}{11/2+\ddots \cfrac{{(-1)}^m mz^2/2}{(2m+1)/2+\ddots}}}}}}} 
\qquad (m\geq 1) \]
Hence the continued fraction expansion given by the theorem follows from 
\[\sum_{n\geq 0} \big| e^{1/l}q_n - p_n \big| \,=\, 2\cdot \frac{\sqrt{\pi}}{2}\cdot \frac{e^{1/l}}{\sqrt{l}}\cdot 
\mbox{erf\,}\big( 1/\sqrt{l} \,\big) \,=\, 2f\big( 1/\sqrt{l} \,\big) \,=\, \frac{2}{l}\, {}_1F_1\big( 1,3/2;1/l \big) \,.\]
\hfill \qed \\
We point out the particular case $z=1$.
\begin{cor}
We have
\begin{eqnarray*}
{}_1F_1\big( 1,3/2;1 \big) &=& e\int_0^1 e^{-t^2}\,dt \\
&=& \mathcal{E}_{MC}(e) \,=\, e-2 + \sum_{n=1}^{\infty} \,\sum_{b=1}^{a_{n}} \big| (bq_{n-1} + q_{n-2})e - (bp_{n-1} + p_{n-2}) \big| \\
&=& \cfrac{1/2}{1/2 - \cfrac{1/2}{3/2 + \cfrac{1}{5/2 - \cfrac{3/2}{7/2+
\cfrac{2}{9/2-\cfrac{5/2}{11/2+\ddots \cfrac{{(-1)}^m m/2}{(2m+1)/2+\ddots}}}}}}} 
\qquad (m\geq 1) \,,
\end{eqnarray*}
where $\mathcal{E}_{MC}(e)$ is the error sum of $e$ taking into account all the minor convergents of 
\begin{equation}
e \,=\, \langle 2; 1,2,1,1,4,1,\dots \rangle \,=\, \langle 2;a_1,a_2,a_3,\dots \rangle \,.
\label{2.1}
\end{equation}
\label{Cor1}
\end{cor}
{\em Proof:\/} \,The formula
\[\mathcal{E}_{MC}(e) \,=\, e\int_0^1 e^{-t^2}\,dt\]
follows by (\ref{2.1}) using
\[\mathcal{E}_{MC}(e) \,=\, e - 1 + \sum_{\nu =0}^{\infty} {(-1)}^{\nu +1} \big( q_{\nu}e - p_{\nu} \big) \Big( \,\frac{1}{2}(1+a_{\nu +1})a_{\nu +1} - a_{\nu +2} \,\Big) \]
and the formulas
\begin{eqnarray*}
q_{3m-1}e - p_{3m-1} &=& -\int_0^1 \frac{x^{m+1}{(x-1)}^m}{m!}e^x\,dx \,,\\
q_{3m}e - p_{3m} &=& -\int_0^1 \frac{x^m{(x-1)}^{m+1}}{m!}e^x\,dx   \,,\\
q_{3m+1}e - p_{3m+1} &=& \int_0^1 \frac{x^{m+1}{(x-1)}^{m+1}}{(m+1)!}e^x\,dx  \\
\end{eqnarray*}
due to H.\,Cohn \cite{Cohn}. \hfill \qed \\
Let  $l\geq 1$. Then we know from \cite[p.\ 193]{Shid} that the numbers $e^{1/l}$ and $\int_0^{1/\sqrt{1/l}}e^{-t^2}\,dt$ are algebraically independent over ${\Q}$. This proves
\begin{cor}
Let $l\geq 2$ be an integer. Then the numbers $\mathcal{E}(e^{1/l}\,)$ and $\mathcal{E}_{MC}(e)$ are transcendental.
\end{cor}

\section{Error sums for $\pi$ and $(\log \rho) /\sqrt{5}$} \label{Sec3}
In \cite{Elsner2}, A.Klauke and the second-named author have found new continued fractions for $1/\pi$ and $(\log \rho)
/\sqrt{5}$. 
In this section we are going to apply these results to compute the corresponding error sums and to decide on their algebraic character. We start with the continued
fraction for $1/\pi$. \\
{\bf 1.)} \,From Theorem\,8 in \cite{Elsner2} and its proof we have the following results.
\begin{equation*}
\begin{aligned}
\frac{1}{\pi}
&=\frac{3}{10} 
\begin{array}{c}\\-\end{array} \frac{14}{25} 
\begin{array}{c}\\-\end{array} \frac{110}{171} 
\begin{array}{c}\\-\ldots -\end{array} \frac{\frac{1}{9}m(m-1)(2m-1)(2m+1)(4m-5)(4m+3)}{ (4m + 1) (4m^2 + 2m-1)} 
\begin{array}{c}\\-\ldots\end{array} \\
&=\frac{p_0}{q_0} 
\begin{array}{c}\\-\end{array} \frac{p_1}{q_1} 
\begin{array}{c}\\-\end{array} \frac{p_2}{q_2} 
\begin{array}{c}\\-\ldots -\end{array} \frac{p_m}{q_m}
\begin{array}{c}\\-\ldots\end{array}  \qquad (m\geq 2)\,.
\end{aligned}
\end{equation*}
Let $n=0,1,2,\dots$. Set
\begin{eqnarray*}
B_n &:=& \frac{2\cdot 4^{n+1}}{n!} \sum_{k=0}^n {n \choose k} (2k+3){\big( k+5/2\big)}_n \,,\\
A_n &:=& \frac{2\cdot 4^{n+1}}{n!} \sum_{k=0}^n \sum_{\nu =0}^k {(-1)}^{k+\nu} {n \choose k} \cfrac{(2k+3){\big( k+
5/2 \big)}_n}{2k-2\nu +1} + {(-4)}^{n+1} \,.
\end{eqnarray*}
Here,
\[{\big( k+5/2 \big)}_n \,=\, (k+5/2)(k+7/2)(k+9/2)\cdots (k+n+3/2) \,.\]
Note that $A_n$ is a rational number, but no integer, while $B_n/4$ is an integer.
Then, for $n\geq 0$, one has
\[\begin{aligned}
\frac{p_0}{q_0} 
\begin{array}{c}\\-\end{array} \frac{p_1}{q_1} 
\begin{array}{c}\\-\end{array} \frac{p_2}{q_2} 
\begin{array}{c}\\-\ldots -\end{array} \frac{p_n}{q_n}
\end{aligned}
\,=\, \frac{B_n}{4A_n} \,,\]
and
\[0\,<\,A_n - \frac{\pi B_n}{4} \,=\, 4(n+3/2)\int_0^1 \frac{t\sqrt{1-t}}{2-t} {\Big( \,\frac{4t(1-t)}{2-t}\,\Big)}^n\,dt \,.\] 
For $0\leq t\leq 1$ the rational function $4t(1-t)/(2-t)$ takes its maximum $2(6-4\sqrt{2})$ at the point $2-\sqrt{2}$. Therefore,
it follows that
\[0\,<\,A_n - \frac{\pi B_n}{4} \,<\, 8(n+3/2){(6-4\sqrt{2})}^n \int_0^1 \frac{t\sqrt{1-t}}{2-t} \,dt 
\,=\, 8(10/3 -\pi)(n+3/2){(6-4\sqrt{2})}^n \,.\] 
The integral on the right-hand side is a Pochhammer integral of a certain hypergeometric function. We show the analogous details
below in part 2.) which is devoted to the number $(\log \rho)/\sqrt{5}$. \\

For (\ref{30}) and (\ref{40}) we define the sequences $B:={(b_n)}_{n\geq 0}$ and $R:={(r_n)}_{n\geq 0}$ by $b_n:=-B_n/4$ and $r_n=-A_n$.
Then we have the error sums
\begin{eqnarray*}
\EE (B,R,\pi) &=& \sum_{m=0}^{\infty} (b_m \pi - r_m) \,=\, -\E (B,R,\pi) \\
&=& -4\int_0^1 \frac{t\sqrt{1-t}}{2-t} \sum_{m=0}^{\infty} (m+3/2) {\Big( \,\frac{4t(1-t)}{2-t}\,\Big)}^m\,dt \\
&=& -4\int_0^1 \frac{t\sqrt{1-t}}{2-t} \cdot \frac{(2-t)(4t^2-7t+6)}{2{(4t^2-5t+2)}^2} \,dt \\
&=& -4\int_0^1 \frac{u^2(1-u^2)(4u^4-u^2+3)}{{(4u^4-3u^2+1)}^2}\,du \,.
\end{eqnarray*}
Here we have introduced the new variable $u:=\sqrt{1-t}$. Computing this integral, we have the following theorem.
\begin{thm}
For the sequences $B:={(b_n)}_{n\geq 0}$ and $R:={(r_n)}_{n\geq 0}$ defined by $b_n:=-B_n/4$ and $r_n=-A_n$ we have
\[\EE (B,R,\pi) \,=\, -\E (B,R,\pi) \,=\, \frac{\sqrt{7}}{49} \log \Big( \,\frac{3-\sqrt{7}}{3+\sqrt{7}}\,\Big) - \frac{3\pi}{2} - \frac{4}{7}
\,=\, - 5.4333111067784 \dots .\] 
\label{Thm2}
\end{thm}
Expressing $\pi$ by $\pi = \frac{2\log i}{i}$, we see that $\E (B,R,\pi)$ is a nonvanishing linear form in logarithms with algebraic arguments
and algebraic coefficients. Then, by Theorem\,2.2 in \cite{Baker}, we have the following corollary.
\begin{cor}
For the sequences $B:={(b_n)}_{n\geq 0}$ and $R:={(r_n)}_{n\geq 0}$ defined by $b_n:=-B_n/4$ and $r_n=-A_n$ the error sum $\E (B,R,\pi)$ is
transcendental, and so is the error sum $\EE (B,R,\pi)$.
\label{Cor2}
\end{cor} 
{\bf 2.)} \,From Theorem\,6 in \cite{Elsner2} and its proof we have the following results.
\begin{equation*}
\begin{aligned}
\frac{\sqrt{5}}{\log \rho}
&=\frac{60}{13} 
\begin{array}{c}\\-\end{array} \frac{7}{80} 
\begin{array}{c}\\-\end{array} \frac{110}{522} 
\begin{array}{c}\\-\ldots -\end{array} \frac{\frac{1}{9}m(m-1)(2m-1)(2m+1)(4m-5)(4m+3)}{ 2(4m + 1) (6m^2 + 3m-1)} 
\begin{array}{c}\\-\ldots\end{array} \\
&=\frac{p_0}{q_0} 
\begin{array}{c}\\-\end{array} \frac{p_1}{q_1} 
\begin{array}{c}\\-\end{array} \frac{p_2}{q_2} 
\begin{array}{c}\\-\ldots -\end{array} \frac{p_m}{q_m}
\begin{array}{c}\\-\ldots\end{array}  \qquad (m\geq 2)\,.
\end{aligned}
\end{equation*}
Let $n=0,1,2,\dots$. Set (cf. (25) in \cite{Elsner2} with $c=d=1$)
\begin{eqnarray*}
D_n &:=& \frac{5\cdot 4^{n+1}}{n!} \sum_{k=0}^n {(-1)}^{n+k}{n \choose k} (2k+3){\big( k+5/2\big)}_n 5^k\,,\\
C_n &:=& 4^n + \frac{4^{n+1}}{n!} \sum_{k=0}^n \sum_{\nu =0}^k {(-1)}^{n+k} {n \choose k} \cfrac{(2k+3){\big( k+
5/2 \big)}_n 5^{\nu}}{2k-2\nu +1} \,.
\end{eqnarray*}
Applying Lemma~6 in \cite{Elsner2} with $x=\tau =1$, we find that
\[\begin{aligned}
\frac{p_0}{q_0} 
\begin{array}{c}\\-\end{array} \frac{p_1}{q_1} 
\begin{array}{c}\\-\end{array} \frac{p_2}{q_2} 
\begin{array}{c}\\-\ldots -\end{array} \frac{p_n}{q_n}
\end{aligned}
\,=\, \frac{D_n}{C_n} \,,\]
and
\[0\,<\,C_n - \frac{D_n\log \rho}{\sqrt{5}} \,=\, \frac{{(5/2)}_n (n+1)!}{4{(5/2)}_{2n+1}} 
\,{}_2F_1 \left( \begin{array}{c} n+1 \quad n+2 \\ 2n+7/2 \end{array} \Big| -\frac{1}{4} \right) \,.\] 
We define the sequences $B:={(b_n)}_{n\geq 0}$ and $R:={(r_n)}_{n\geq 0}$ by $b_n:=-D_n$ and $r_n=-C_n$.
Then we have the error sum
\[\EE \Big( \,B,R,\frac{\log \rho}{\sqrt{5}}\,\Big) \,=\, \sum_{m=0}^{\infty} \Big( \,b_m \frac{\log \rho}{\sqrt{5}} - r_m\,\Big) \\
= -\frac{1}{4} \sum_{m=0}^{\infty} \frac{{(5/2)}_m (m+1)!}{4{(5/2)}_{2m+1}} 
\,{}_2F_1 \left( \begin{array}{c} m+1 \quad m+2 \\ 2m+7/2 \end{array} \Big| -\frac{1}{4} \right) \,.\]
To compute this error sum, the method is the same as used above for the error sum of $\E(\pi)$. First we express the hypergeometric function by Pochhammer's 
integral. Let $a,b,c,z$ be complex numbers satisfying $|z|<1$, $\Re(c-b)>0$, and $\Re(b)>0$. Then we have the identity
\[{}_2F_1 \left( \begin{array}{c} a \quad b \\ c \end{array} \Big| z \right) \,=\, \frac{\Gamma(c)}{\Gamma(b) \Gamma(c-b)}
\int_0^1 t^{b-1} {(1-t)}^{c-b-1} {(1-zt)}^{-a}\,dt \,,\]
cf.\,\cite[p.\ 20]{Slater}. The conditions are fulfilled for $a=n+1$, $b=n+2$, $c=2n+7/2$, and $z=-1/4$, where $n=0,1,2,\dots$.
Hence, it follows that
\[{}_2F_1 \left( \begin{array}{c} n+1 \quad n+2 \\ 2n+7/2 \end{array} \Big| -\frac{1}{4} \right) \,=\, \frac{\Gamma(2n+7/2)}
{\Gamma(n+2) \Gamma(n+3/2)} \int_0^1 t^{n+1}{(1-t)}^{n+1/2}{\Big( \,1+\frac{t}{4}\,\Big)}^{-n-1} \,dt \,.\]
In order to simplify the above expressions we need two identities involving Pochhammer's symbol (\cite[p.\ 239]{Slater}).
\begin{eqnarray*}
\frac{{(5/2)}_n}{{(5/2)}_{2n+1}} &=& \frac{1}{{(n+5/2)}_{n+1}} \,,\\
\frac{\Gamma(2n+7/2)}{\Gamma(n+3/2)} &=& (n+3/2){(n+5/2)}_{n+1} \,. 
\end{eqnarray*}
Collecting together all the above results, it follows that
\begin{eqnarray*}
\EE \Big( \,B,R,\frac{\log \rho}{\sqrt{5}}\,\Big) &=& -\frac{1}{4} \sum_{m=0}^{\infty} (m+3/2) 
\int_0^1 \frac{t^{m+1}{(1-t)}^{m+1/2}}{{(1+t/4)}^{m+1}} \,dt \\
&=& -\frac{1}{4} \int_0^1 \frac{4t\sqrt{1-t}}{4+t} \sum_{m=0}^{\infty} (m+3/2) {\Big( \,\frac{4t(1-t)}{4+t}\,\Big)}^m \,dt \\
&=& -\frac{1}{2} \int_0^1 \frac{t(4t^2-t+12)\sqrt{1-t}}{{(4t^2-3t+4)}^2} \,dt \\
&=& \int_0^1 \frac{u^2(1-u^2)(4u^4-7u^2+15)}{{(4u^4-5u^2+5)}^2} \,du \,, 
\end{eqnarray*}
where $u=\sqrt{1-t}$. This proves
\begin{thm}
For the sequences $B:={(b_n)}_{n\geq 0}$ and $R:={(r_n)}_{n\geq 0}$ defined by $b_n:=-D_n$ and $r_n=-C_n$ we have
\begin{eqnarray*}
&& \EE \Big( \,B,R,\frac{\log \rho}{\sqrt{5}}\,\Big) \\
&=& \cfrac{\sqrt{124\sqrt{5}-265}\log \Big( \,1+\cfrac{\sqrt{5}}{2} - \cfrac{\sqrt{4\sqrt{5}+5}}{2}\,\Big) - \sqrt{124\sqrt{5}+265}\arccos \Big( \,\cfrac{\sqrt{5}}{2}-1\,\Big)}
{55\sqrt{11}} + \frac{1}{11} \\ \\
&=& - 0.1210649459927\dots \,.
\end{eqnarray*}
\label{Thm3}
\end{thm}
Using $\arccos z = \frac{1}{i} \log (z+\sqrt{z^2-1})$, we obtain by Theorem\,2.2 in \cite{Baker} the following result.
\begin{cor}
For the sequences $B:={(b_n)}_{n\geq 0}$ and $R:={(r_n)}_{n\geq 0}$ defined by $b_n:=-D_n$ and $r_n=-C_n$ the error sum $\E (B,R,\log \rho /\sqrt{5})$ is
transcendental.
\label{Cor3}
\end{cor}

\section{An error sum for $\log (1+t)$} \label{Sec4}
In this section we generalize a concept from the  proof of Theorem~3 in \cite{Elsner3}, where a nonregular continued 
fraction for $\log 2$ is established. First we shall prove a continued fraction expansion for $\log (1+t)$ with $-1<t\leq 1$, namely
\begin{equation}
\log (1+t) \,=\, \frac{2t}{2+t} 
\begin{array}{c}\\-\end{array} \frac{1^2t^2}{3(2+t)} 
\begin{array}{c}\\-\end{array} \frac{2^2t^2}{5(2+t)}
\begin{array}{c}\\-\end{array} \frac{3^2t^2}{7(2+t)}
\begin{array}{c}\\- \dots -\end{array} \frac{m^2t^2}{(2m+1)(2+t)} 
\begin{array}{c}\\- \dots \end{array} \,,
\label{4.10}
\end{equation}
where $m=1,2,\dots$. O.Perron \cite[p.\ 152]{Perron} cites by equation (7) the continued fraction
\[\log (1+t) \,=\, \frac{t}{1} 
\begin{array}{c}\\+\end{array} \frac{1^2t}{2} 
\begin{array}{c}\\+\end{array} \frac{1^2t}{3}
\begin{array}{c}\\+\end{array} \frac{2^2t}{4}
\begin{array}{c}\\+\end{array} \frac{2^2t}{5}
\begin{array}{c}\\+\end{array} \frac{3^2t}{6} 
\begin{array}{c}\\+\end{array} \frac{3^2t}{7}  
\begin{array}{c}\\+ \dots \end{array}  
\,.\]
Here we shall give full details of the proof, since a new argument is needed in H.Cohen's method  \cite{Cohen} established for Ap{\'e}ry's
irregular continued fractions of $\zeta(2)$ and $\zeta(3)$, and we need the details in order to compute the error sum. Similar to Ap{\'e}ry's approach we have to handle with 
combinatorial series.  \\
In the sequel we fix a real number $t$ with  $-1<t\leq 1$. Let $n\geq 0$ be an integer. We define two combinatorial series by
\begin{eqnarray*}
B_n &:=& \sum_{k=0}^n {n \choose k}{n+k \choose k} t^{n-k} \,,\\
A_n &:=& \sum_{k=0}^n {n \choose k}{n+k \choose k} t^{n-k} c_k  \,,
\end{eqnarray*}
where 
\[c_k \,:=\, \sum_{m=1}^k \frac{{(-1)}^{m-1}t^m}{m} \,.\]
By applying {\em Zeilberger's algorithm\/} \cite[Ch.\ 7]{Koepf} (Algorithm\,7.1) on a computer algebra system, it turns out that the numbers $B_n$ satisfy the linear
three-term recurrence formula
\begin{equation}
 (n+1)X_{n+1} -(2n+1) (2+t)X_n + nt^2X_{n-1} \,=\, 0 \qquad (n\geq 1)\,.
\label{4.20}
\end{equation}
In the sequel we prove this formula for $X_n=B_n$ without using a computer, since we need the details to show that even
$X_n=A_n$ satisfy the recurrence. Let $k,n$ denote integers. Set
\begin{eqnarray*}
\lambda_{n,k} &:=& {n \choose k}{n+k \choose k} t^{n-k} \,,\\
B_{n,k} &:=& -(4n+2)\lambda_{n,k} \,,\\
A_{n,k} &:=& B_{n,k}c_k \,,\\
S_{n,k} &:=& (n+1)\lambda_{n+1,k}c_k - (2n+1)(2+t)\lambda_{n,k}c_k + nt^2\lambda_{n-1,k}c_k \,.
\end{eqnarray*}
Note that ${n \choose k}=0$ for $k<0$ or $k>n$, which implies that $A_{n,n+1}=B_{n,n+1}=A_{n,-1}=B_{n,-1}=0$. One easily verfies the identities\footnote{We should like to point 
out that there is a misprint in the formula for $\lambda_{n+1,k}/\lambda_{n,k}$ in \cite{Elsner3}.}
\begin{eqnarray*}
\frac{\lambda_{n,k-1}}{\lambda_{n,k}} &=& \frac{k^2t}{(n+k)(n-k+1)} \,,\\
\frac{\lambda_{n+1,k}}{\lambda_{n,k}} &=& \frac{(n+k+1)t}{n-k+1} \,,\\
\frac{\lambda_{n-1,k}}{\lambda_{n,k}} &=& \frac{n-k}{(n+k)t} \,,
\end{eqnarray*}
which can be applied to prove the identity
\begin{equation}
B_{n,k} - B_{n,k-1} \,=\, (n+1)\lambda_{n+1,k} - (2n+1)(2+t)\lambda_{n,k} + nt^2\lambda_{n-1,k} \,.
\label{4.30}
\end{equation}
Summing up on both sides of (\ref{4.30}) from $k=0$ to $k=n+1$, we obtain
\[0 \,=\, B_{n,n+1} - B_{n,-1} \,=\, \sum_{k=0}^{n+1} \big( B_{n,k} - B_{n,k-1}  \big) \,=\, (n+1)B_{n+1} - (2n+1)(2+t)B_n + nt^2B_{n-1} \,,\]
which proves (\ref{4.20}) for $X_n=B_n$. \\
Multiplying (\ref{4.30}) by $c_k$, we obtain $S_{n,k}=(B_{n,k}-B_{n,k-1})c_k$. Hence,
\begin{eqnarray*}
A_{n,k} - A_{n,k-1} &=& B_{n,k}c_k - B_{n,k-1}c_{k-1} \,=\, \big( B_{n,k} - B_{n,k-1} \big) c_k + B_{n,k-1} \big( c_k - c_{k-1} \big) \\
&=& S_{n,k} + B_{n,k-1} \frac{{(-1)}^{k-1}t^k}{k} \,.
\end{eqnarray*}
Again, we sum up from $k=0$ to $k=n+1$. This gives
\begin{eqnarray*}
0 &=& A_{n,n+1} - A_{n,-1} \,=\, \sum_{k=0}^{n+1} \big( A_{n,k} - A_{n,k-1} \big) \\
&=& \sum_{k=0}^{n+1} S_{n,k} - (4n+2)\sum_{k=1}^{n+1} {n \choose k-1}{n+k-1 \choose k-1} t^{n-k+1}
\frac{{(-1)}^{k-1}t^k}{k} \\
&=& (n+1)A_{n+1} - (2n+1)(2+t)A_n + nt^2A_{n-1} - (4n+2)t^{n+1}\sum_{k=0}^n {n \choose k}{n+k \choose k}
\frac{{(-1)}^k}{k+1} \,.
\end{eqnarray*}
Finally Vandermonde's theorem for the hypergeometric series ${}_2F_1(n+1,-n,2;1)$ (\cite[eq.\ (1.7.7)]{Slater}) completes our proof of (\ref{4.20}) for $X_n=A_n$  
by
\[\sum_{k=0}^n {n \choose k}{n+k \choose k} \frac{{(-1)}^k}{k+1} \,=\, {}_2F_1(n+1,-n,2;1) \,=\, \frac{{(1-n)}_n}{{(2)}_n} \,=\, 0 \]
for $n\geq 1$. In the next step we prove that
\begin{equation}
\lim_{n\to \infty} \frac{A_n}{B_n} \,=\, \log (1+t) \,.
\label{4.40}
\end{equation}
For this purpose we shall prove that for every fixed integer $\nu \geq 0$ we have the limit
\begin{equation}
\lim_{n\to \infty} \frac{\displaystyle{n \choose \nu}{n+\nu \choose \nu}t^{n-\nu}}{\displaystyle \sum_{k=0}^n 
{n \choose k}{n+k \choose k} t^{n-k}} \,=\, 0 \,.
\label{4.50}
\end{equation}
Then, (\ref{4.50}) implies (\ref{4.40}) by a theorem of O.Toeplitz (\cite[p.\ 10, no.\ 66]{Polya}), since
\[\lim_{n\to \infty} c_n \,=\, \sum_{m=1}^{\infty} \frac{{(-1)}^{m-1}t^m}{m} \,=\, \log (1+t) \,.\]
There is nothing to show for $t=0$, because $A_n=0$ and $B_n={2n \choose n}\not= 0$. Therefore, keep $\nu \in {\N}_0$
and $t\in (-1,1] \setminus \{0\}$ fixed. We substitute $X_n=B_n$ into (\ref{4.20}) and divide the equation by $(n+1)B_n$.
Then we obtain
\[\frac{B_{n+1}}{B_n} - \frac{(2n+1)(2+t)}{n+1} + \frac{nt^2}{n+1}\cdot \frac{1}{B_n/B_{n-1}} \,=\, 0 \,.\]
Let $\alpha :=\lim_{n\to \infty} B_{n+1}/B_n$. By taking the limit $n\to \infty$, it follows that $\alpha$ satisfies the quadratic
equation
\[\alpha - 2(2+t) + \frac{t^2}{\alpha} \,=\, 0\,,\]
which yields
\[\alpha \,=\, 2+t+2\sqrt{1+t} \,>\, |t| \qquad (-1<t\leq 1) \,.\]
Put $\beta :=(\alpha +|t|)/2$. Then,
\begin{equation}
0 \,<\,|t| \,<\, \beta \,<\, \alpha \,.
\label{4.60}
\end{equation}
There is an integer $n_0=n_0(t)$ satisfying
\[\frac{B_m}{B_{m-1}} \,>\, \beta \,>\, 0 \qquad (m\geq n_0) \,,\]
or
\[|B_m| \,>\, \beta |B_{m-1}| \qquad (m\geq n_0) \,.\]
Then, for $n\geq 2n_0$ and $k:=n-n_0\geq n_0$, we have
\[|B_n| \,>\, \beta |B_{n-1}| \,>\, \beta^2 |B_{n-2}| \,> \dots >\beta^k |B_{n-k}| \,=\, \beta^{n-n_0} |B_{n_0}| \,.\]
Consequently, we obtain
\begin{eqnarray*}
&& \lim_{n\to \infty}\left| \,\frac{\displaystyle{n \choose \nu}{n+\nu \choose \nu}t^{n-\nu}}{\displaystyle \sum_{k=0}^n 
{n \choose k}{n+k \choose k} t^{n-k}} \,\right| \,=\, \lim_{n\to \infty} \frac{\displaystyle{n \choose \nu}{n+\nu \choose \nu}
{|t|}^{n-\nu}}{|B_n|} \\ \\
&\leq &   \lim_{n\to \infty} \frac{\displaystyle{n \choose \nu}{n+\nu \choose \nu}
{|t|}^{n-\nu}}{\beta^{n-n_0}|B_{n_0}|} \,=\, \frac{\beta^{n_0}}{{|t|}^{\nu}|B_{n_0}|} \lim_{n\to \infty}
{n \choose \nu}{n+\nu \choose \nu}{\Big( \,\frac{|t|}{\beta}\,\Big)}^n  \\  \\
&=& 0 \,,
\end{eqnarray*}
since $0<|t|{\beta}^{-1} <1$ by (\ref{4.60}), and ${n \choose \nu}{n+\nu \choose \nu}$ is a polynomial in $n$ of degree
$2\nu$. This completes the proof of (\ref{4.50}) and, consequently, of (\ref{4.40}). \\
We rewrite the recurrence formula (\ref{4.20}) as
\[P(n+1)X_{n+1} - Q(n+1)X_n - R(n+1)X_{n-1} \,=\, 0\,,\]
where
\begin{eqnarray*}
P(n+1) &:=& n+1 \,,\\
Q(n+1) &:=& (2n+1)(2+t) \,,\\
R(n+1) &:=& -nt^2 \,.
\end{eqnarray*}
Then, we obtain
\[\log (1+t) \,=\, b_0 + \frac{a_1}{b_1} 
\begin{array}{c}\\+\end{array} \frac{a_2}{b_2}
\begin{array}{c}\\+\end{array} \frac{a_3}{b_3}
\begin{array}{c}\\+ \dots  \end{array} \]
with
\[\begin{array}{lcrlrllc}
b_0 &=& 0\,,\quad b_1 &=& 2+t \,,\quad b_{n+1} &=& \displaystyle \frac{Q(n+1)}{P(n+1)} \,=\, \frac{(2n+1)(2+t)}{n+1} \,,\\ \\
&& a_1 &=& 2t \,,\quad a_{n+1} &=& \displaystyle \frac{R(n+1)}{P(n+1)} \,=\, -\frac{nt^2}{n+1} \,.
\end{array} \]
This gives the continued fraction
\[\log(1+t) \,=\, \frac{2t}{2+t} 
\begin{array}{c}\\-\end{array} \frac{t^2/2}{3(2+t)/2}
\begin{array}{c}\\-\end{array} \frac{2t^2/3}{5(2+t)/3}
\begin{array}{c}\\-\end{array} \frac{3t^2/4}{7(2+t)/4}
\begin{array}{c}\\- \dots \end{array} \,,\]
which is equivalent with (\ref{4.10}). \\
Next, we compute the error sum $\EE (B,R,\log (1+t))=\sum_{m=0}^{\infty} (B_m\log (1+t) -A_m)$ for $B:={(B_n)}_{n\geq 0}$ and 
$R:={(A_n)}_{n\geq 0}$.
\begin{lemma}
Let $-1<t\leq 1$. For every integer $n\geq 0$ we have
\[B_n\log (1+t) - A_n \,=\, t^{2n+1} \int_0^1 \frac{x^n{(1-x)}^n}{{(1+tx)}^{n+1}} \,dx \,.\] 
\label{Lem4.1}
\end{lemma}
{\em Proof.\/} \,For $A_n$ and $B_n$ defined above, we obtain
\begin{eqnarray*}
&& B_n\log (1+t) - A_n \\
&=& \sum_{k=0}^n {n \choose k}{n+k \choose k} t^{n-k} \sum_{m=k+1}^{\infty} \frac{{(-1)}^{m-1}t^m}{m} \\
&=& \sum_{m=0}^{\infty} {(-1)}^m t^{m+1} \sum_{k=0}^n {(-1)}^k {n \choose k}{n+k \choose k} \frac{t^{n-k}t^k}{m+k+1} \\ 
&=& t^{n+1} \int_0^1 \sum_{m=0}^{\infty} {(-1)}^m {(tx)}^m \sum_{k=0}^n {(-1)}^k {n \choose k}{n+k \choose k} x^k\,dx \\
&=& t^{n+1} \int_0^1 \sum_{m=0}^{\infty} {(-tx)}^m \frac{d^n}{dx^n} \Big( \,\frac{x^n{(1-x)}^n}{n!}\,\Big) \,dx \\
&=& t^{n+1} \int_0^1 \frac{1}{1+tx} \cdot \frac{d^n}{dx^n} \Big( \,\frac{x^n{(1-x)}^n}{n!}\,\Big) \,dx \\
&=& {(-1)}^n t^{n+1} \int_0^1 \frac{d^n}{dx^n} \Big( \,\frac{1}{1+tx}\,\Big) \cdot \frac{x^n{(1-x)}^n}{n!} \,dx \\
&=& t^{2n+1} \int_0^1 \frac{x^n{(1-x)}^n}{{(1+tx)}^{n+1}} \,dx \,.
\end{eqnarray*}
The antiderivative in the last but one line was obtained using $n$-fold integration by parts. The lemma is proven. \hfill \qed

A first consequence of Lemma~\ref{Lem4.1} is an explicit formula for the error sum of $\log (1+t)$. 
\begin{cor}
Let $-1<t\leq 1$. For the sequences $B:={(B_n)}_{n\geq 0}$ and $R:={(A_n)}_{n\geq 0}$ we have 
\[\EE \big( B,R,\log (1+t)\big) \,=\, \frac{2}{\sqrt{3+2t-t^2}} \left( \arctan \Big( \,\frac{1+t}{\sqrt{3+2t-t^2}} - 
\arctan \Big( \,\frac{1-t}{\sqrt{3+2t-t^2}}\,\Big) \right) \,.\]
In particular, $\EE (B,R,\log 2)=\pi/4$.
\label{Cor4.1}
\end{cor}
{\em Proof:\/} \,From Lemma~\ref{Lem4.1} we obtain
\begin{eqnarray*}
&& \EE \big( B,R,\log (1+t)\big) \\
&=& t\int_0^1 \sum_{m=0}^{\infty} \frac{t^{2m}x^m{(1-x)}^m}{{(1+tx)}^{m+1}} \,dx \\
&=& t\int_0^1 \frac{1}{1+tx} \sum_{m=0}^{\infty} {\Big( \,\frac{t^2x(1-x)}{1+tx}\,\Big)}^m\,dx \\
&=& t\int_0^1 \frac{dx}{1+t(1-t)x+t^2x^2} \\
&=& \frac{2}{\sqrt{3+2t-t^2}} \left( \arctan \Big( \,\frac{1+t}{\sqrt{3+2t-t^2}} - \arctan \Big( \,\frac{1-t}{\sqrt{3+2t-t^2}}\,\Big) \right) \,.
\end{eqnarray*}
This proves the corollary. \hfill \qed
\\
By straightforward computations it can be seen that the error sum from Corollary~\ref{Cor4.1} satisfies a linear first order differential equation.
\begin{cor}
Let $-1<t\leq 1$. For the sequences $B:={(B_n)}_{n\geq 0}$ and $R:={(A_n)}_{n\geq 0}$ the function $f(t):=\EE (B,R,\log (1+t))$ satisfies the differential equation
\[(3+2t-t^2)f'+(1-t)f-3 \,=\, 0 \,,\]
where $f'=df/dt$.
\label{Cor4.2}
\end{cor}
A second consequence of Lemma~\ref{Lem4.1} is 
\[\EE \big( B,R,\log (1+t)\big) \,=\, \mbox{sign}(t) \E \big( B,R,\log (1+t)\big) \,.\]
Finally, the continued fraction (\ref{4.10}) and Lemma~\ref{Lem4.1} allow to prove the irrationality of $\log (1+t)$ for certain rationals $t:=a/b$. 
\begin{cor}
Let $0<a/b\leq 1$ be a rational number with $ea^2<4b$. Then the number $\log (1+a/b)$ is irrational. In particular, for every integer $k\geq 1$ the number $\log (1+1/k)$
is irrational.
\label{Cor4.3}
\end{cor} 
{\em Proof:\/} \,Let $d_n :=\mbox{l.c.m.}(1,2,3,\dots,n)$ denote the least common multiple of the integers $1,2,3,\dots ,n$. One
knows by the prime number theorem that
\[\log d_n \,=\, \sum_{p\leq n} \Big[ \,\frac{\log n}{\log p}\,\Big] \log p \,\sim \, n \,,\] 
where $p$ runs through all primes less than or equal to $n$ (\cite[Theorem\,434]{Hardy}). By the hypothesis $ea^2<4b$ there is a positive real number $\varepsilon$ 
such that $e^{1+\varepsilon}a^2<4^{1-\varepsilon}b$. Hence, for all sufficiently large numbers $n$, it follows that
\[d_n \,<\, e^{(1+\varepsilon)n} \,.\]
Let $t=a/b$. With $b^nd_nA_n \in {\Z}$ and $b^nd_nB_n \in {\Z}$ we know by Lemma~\ref{Lem4.1} that
\begin{eqnarray*}
0 &<& \big| b^nd_nB_n\log (1+t) - b^nd_nA_n \big| \\
&=& \frac{a}{b} b^nd_n {\Big( \,\frac{a}{b}\,\Big)}^{2n} \int_0^1 \frac{x^n{(1-x)}^n}{{(1+ax/b)}^{n+1}} \,dx\\
&<& \frac{a}{b} \cdot \frac{e^{(1+\varepsilon)n}a^{2n}}{b^n} \int_0^1 x^n{(1-x)}^n \,dx \\
&<& \frac{a}{b} \cdot 4^{(1-\varepsilon)n} \int_0^1 \frac{dx}{4^n}\,dx \\
&=& \frac{t}{4^{\varepsilon n}} \,\to \, 0 
\end{eqnarray*}
for $n$ tending to infinity. This completes the proof of Corollary~\ref{Cor4.3}. \hfill \qed

\section{On error sums formed by Ap{\'e}ry's continued fractions for $\zeta(2)$ and $\zeta(3)$} \label{Sec5}
Computing the error sums formed by the linear three term recurrences and continued fractions of $\zeta(2)$, $\zeta(3)$ introduced by R.\,Ap{\'e}ry, this
leads unexpectedly into a wide field of connections between famous sequences of integers. For the needed results we refer to \cite{Cohen} and \cite{Beukers}. \\
{\bf 1.)} Error sums for $\zeta(2)$. We have \\
\begin{eqnarray*}
\zeta(2) &=& \frac{\pi^2}{6} \,=\, \frac{5}{3} 
\begin{array}{c}\\+\end{array} \frac{1^4}{25}
\begin{array}{c}\\+\end{array} \frac{2^4}{69}
\begin{array}{c}\\+ \dots +\end{array} \frac{n^4}{11n^2+11n+3}
\begin{array}{c}\\+ \dots \end{array} \\
&=& b_0 + \frac{a_1}{b_1} 
\begin{array}{c}\\+\end{array} \frac{a_2}{b_2}
\begin{array}{c}\\+\end{array} \frac{a_3}{b_3}
\begin{array}{c}\\+ \dots \end{array}
\end{eqnarray*}
with
\[b_0=0 \,,\quad b_1=3 \,,\quad b_{n+1} = 11n^2+11n+3 \quad (n\geq 1) \,,\]
\[a_1=5 \,,\quad a_{n+1} = n^4 \quad (n\geq 1) \,.\]
A recurrence formula for both sequences
\begin{eqnarray*}
B_n &:=& \sum_{k=0}^n {n \choose k}^2 {n+k \choose k} \,,\\
A_n &:=& \sum_{k=0}^n {n \choose k}^2 {n+k \choose k}\left( 2\sum_{m=1}^n \frac{{(-1)}^{m-1}}{m^2} + \sum_{m=1}^k \cfrac{{(-1)}^{n+m-1}}{m^2{n \choose m}{n+m \choose m}} \right) \,,
\end{eqnarray*}
is
\[0 \,=\, {(n+1)}^2X_{n+1} - (11n^2+11n+3)X_n - n^2X_{n-1} \,.\]
Then,
\[b_0 + \frac{a_1}{b_1} 
\begin{array}{c}\\+\end{array} \frac{a_2}{b_2}
\begin{array}{c}\\+\end{array} \frac{a_3}{b_3}
\begin{array}{c}\\+ \dots +\end{array} \frac{a_n}{b_n} \\
\,=\, \frac{A_n}{B_n} \,.\]
We obtain from \cite[eq.\ (5)]{Beukers} for the sequences $B_2:={(B_n)}_{n\geq 0}$ and $R_2:={(A_n)}_{n\geq 0}$,
\begin{eqnarray}
\E^* \big( B_2,R_2,\zeta(2)\big) &=& \sum_{n=0}^{\infty} \Big( \,B_n\zeta(2) - A_n\, \Big) \,=\, \sum_{n=0}^{\infty} {(-1)}^n \int_0^1 \int_0^1 \frac{x^n{(1-x)}^ny^n{(1-y)}^n}{{(1-xy)}^{n+1}} \,dx\,dy
\nonumber \\
&=& \int_0^1 \int_0^1 \frac{dx\,dy}{1+x^2y^2-xy^2-yx^2} \,=\,  1.5832522167 \dots \,.
\label{5.10} 
\end{eqnarray}
Similarly, one has
\begin{eqnarray}
\E\big( B_2,R_2,\zeta(2)\big)  &=& \sum_{n=0}^{\infty} \Big| \,B_n\zeta(2) - A_n\, \Big| \,=\, \sum_{n=0}^{\infty} \int_0^1 \int_0^1 \frac{x^n{(1-x)}^ny^n{(1-y)}^n}{{(1-xy)}^{n+1}} \,dx\,dy \nonumber 
\\
&=& \int_0^1 \int_0^1 \frac{dx\,dy}{1-x^2y^2-2xy+xy^2+yx^2} \,=\, 1.7141459142 \dots \,. 
\label{5.20}
\end{eqnarray}
{\bf 2.)} Error sums for $\zeta(3)$. Here,  \\
\begin{eqnarray*}
\zeta(3) &=& \frac{6}{5} 
\begin{array}{c}\\-\end{array} \frac{1^6}{117}
\begin{array}{c}\\-\end{array} \frac{2^6}{535}
\begin{array}{c}\\- \dots -\end{array} \frac{n^6}{34n^3+51n^2+27n+5}
\begin{array}{c}\\- \dots \end{array} \\
&=& b_0 + \frac{a_1}{b_1} 
\begin{array}{c}\\+\end{array} \frac{a_2}{b_2}
\begin{array}{c}\\+\end{array} \frac{a_3}{b_3}
\begin{array}{c}\\+ \dots \end{array}
\end{eqnarray*}
with
\[b_0=0 \,,\quad b_1=5 \,,\quad b_{n+1} = 34n^3+51n^2+27n+5 \quad (n\geq 1) \,,\]
\[a_1=6 \,,\quad a_{n+1} = -n^6 \quad (n\geq 1) \,.\]
A recurrence formula for both sequences, 
\begin{eqnarray*}
D_n &:=& \sum_{k=0}^n {n \choose k}^2 {n+k \choose k}^2 \,,\\
C_n &:=& \sum_{k=0}^n {n \choose k}^2 {n+k \choose k}^2\left( \sum_{m=1}^n \frac{1}{m^3} + \sum_{m=1}^k \cfrac{{(-1)}^{m-1}}{2m^3{n \choose m}{n+m \choose m}} \right) \,,
\end{eqnarray*}
is
\begin{eqnarray*}
0 &=& P(n+1)X_{n+1} - Q(n+1)X_n - R(n+1)X_{n-1} \\
&=& {(n+1)}^3X_{n+1} - (34n^3+51n^2+27n+5)X_n + n^3X_{n-1} \,.
\end{eqnarray*}
The construction of $C_n$ and $D_n$ leads to the identity
\[b_0 + \frac{a_1}{b_1} 
\begin{array}{c}\\+\end{array} \frac{a_2}{b_2}
\begin{array}{c}\\+\end{array} \frac{a_3}{b_3}
\begin{array}{c}\\+ \dots +\end{array} \frac{a_n}{b_n} \\
\,=\, \frac{C_n}{D_n} \,.\]
We obtain from \cite[eq.\ (7)]{Beukers} for the sequences $B_3:={(D_n)}_{n\geq 0}$ and $R_3:={(C_n)}_{n\geq 0}$,
\begin{eqnarray}
&& \E^*\big( B_3,R_3,\zeta(3)\big) \nonumber \\
&=& \sum_{n=0}^{\infty} \Big( \,D_n\zeta(3) - C_n\, \Big) \,=\, \sum_{n=0}^{\infty} \Big| \,D_n\zeta(3) - C_n\, \Big| \,=\, \E\big( B_3,R_3,\zeta(3)\big) \nonumber \\
&=& \sum_{n=0}^{\infty} \frac{1}{2} \int_0^1 \int_0^1 \int_0^1 \frac{x^n{(1-x)}^ny^n{(1-y)}^nw^n{(1-w)}^n}{{\big( 1-(1-xy)w\big)}^{n+1}} \,dx\,dy\,dw \nonumber \\
&=& \frac{1}{2} \int_0^1 \int_0^1 \int_0^1 \frac{dx\,dy\,dw}{1+x^2y^2w^2-xy^2w^2-x^2yw^2-x^2y^2w+xyw^2+xy^2w+x^2yw-w} 
\nonumber \\ \nonumber \\
&=& 1.2124982529 \dots \,. 
\label{5.30}
\end{eqnarray}
{\bf 3.)} Now we focus our interest on various methods in order to express the multiple integrals in (\ref{5.10}), (\ref{5.20}), and (\ref{5.30}), by series with rational terms. 
A first approach to this subject involves the hypergeometric function. 
\begin{thm}
For the sequences $B_i,R_i$ $(i=2,3)$ defined above for $\zeta(2)$ and $\zeta(3)$, respectively, we have
\begin{eqnarray*}
\E\big( B_2,R_2,\zeta(2)\big) &=& \sum_{n=0}^{\infty} \,\,\sum_{k=0}^{\infty} \cfrac{\displaystyle {n+k \choose n}}{\displaystyle (2n+k+1)^2{2n+k \choose n}^2} \\
&=& \sum_{n=0}^{\infty} \frac{\displaystyle {}_3F_2\left( \begin{array}{c} n+1 \quad n+1 \quad n+1 \\ 2n+2 \quad 2n+2  \end{array} \,\Big| 1 \right) }{\displaystyle{(2n+1)}^2{2n \choose n}^2} 
\,,\\
\E^*\big( B_2,R_2,\zeta(2)\big) &=& \sum_{n=0}^{\infty} \,\,\sum_{k=0}^{\infty} \cfrac{\displaystyle {(-1)}^n  {n+k \choose n}}{\displaystyle (2n+k+1)^2{2n+k \choose n}^2} \\
&=& \sum_{n=0}^{\infty} {(-1)}^n \frac{\displaystyle {}_3F_2\left( \begin{array}{c} n+1 \quad n+1 \quad n+1 \\ 2n+2 \quad 2n+2  \end{array} \,\Big| 1 \right)}
{\displaystyle{(2n+1)}^2{2n \choose n}^2} \,,\\
\E\big( B_3,R_3,\zeta(3)\big) &=& \frac{1}{2}\sum_{n=0}^{\infty} \,\,\sum_{k=0}^{\infty} \,\,\sum_{l=0}^k \cfrac{\displaystyle {(-1)}^l  {k \choose l}{n+l \choose n}}
{\displaystyle (2n+k+1)(2n+l+1)^2{2n+k \choose n}{2n+l \choose n}^2} \\
&=&  \frac{1}{2}\sum_{n=0}^{\infty}\,\,\sum_{k=0}^{\infty} \frac{\displaystyle {}_4F_3\left( \begin{array}{c} n+1 \quad n+1 \quad n+1 \quad -k \\ 2n+2 \quad 2n+2 \quad 1  \end{array} \,\Big| 
1 \right)}{\displaystyle{(2n+1)}^2(2n+k+1){2n \choose n}^2{2n+k \choose n}} \,.
\end{eqnarray*}
\label{Thm5.10}
\end{thm}
Note that the hypergeometric function
\[{}_4F_3\left(\begin{array}{c} n+1 \quad n+1 \quad n+1 \quad -k \\ 2n+2 \quad 2n+2 \quad 1  \end{array} \,\Big| 1 \right) \]
takes rational values for all $0\leq k,n<\infty$. \\
{\em Proof:\/} \,It suffices to prove the identities for $\E\big( B_2,R_2,\zeta(2)\big)$, since the arguments are the same for the remaining error sums. 
The basic idea is to use the expansion
\[\frac{1}{{(1-t)}^{n+1}} \,=\, \sum_{k=0}^{\infty} \frac{{(n+1)}_k}{k!} t^k \,.\]
Then, (\ref{5.20}) gives
\begin{eqnarray}
\E\big( B_2,R_2,\zeta(2)\big) &=& \sum_{n=0}^{\infty} \int_0^1 \int_0^1 \frac{x^n{(1-x)}^n y^n{(1-y)}^n}{{(1-xy)}^{n+1}}\,dx\,dy \nonumber \\
&=& \sum_{n=0}^{\infty} \int_0^1 \int_0^1 \sum_{k=0}^{\infty} \frac{{(n+1)}_k}{k!} {(xy)}^k x^n{(1-x)}^n y^n{(1-y)}^n \,dx\,dy \nonumber \\
&=& \sum_{n=0}^{\infty} \,\,\sum_{k=0}^{\infty} \frac{{(n+1)}_k}{k!} \int_0^1 x^{n+k}{(1-x)}^n\,dx \int_0^1 y^{n+k}{(1-y)}^n\,dy \nonumber \\
&=& \sum_{n=0}^{\infty} \,\,\sum_{k=0}^{\infty} \frac{{(n+1)}_k}{k!} {\Big( \,\frac{\Gamma(n+1) \Gamma(n+k+1)}{\Gamma(2n+k+2)}\, \Big)}^2 
\label{5.40} \\
&=& \sum_{n=0}^{\infty} \,\,\sum_{k=0}^{\infty} \frac{n! {(n+k)!}^3}{k! {(2n+k)!}^2{(2n+k+1)}^2} 
\,=\, \sum_{n=0}^{\infty} \,\,\sum_{k=0}^{\infty} \cfrac{\displaystyle {n+k \choose n}}{\displaystyle (2n+k+1)^2{2n+k \choose n}^2} \,.
\nonumber 
\end{eqnarray}
The second identity for $\E\big( B_2,R_2,\zeta(2)\big)$ in Theorem~\ref{Thm5.10} follows from (\ref{5.40}) and from 
\[\frac{\Gamma^2(n+1) \Gamma^2(n+k+1)}{\Gamma^2(2n+k+2)} \,=\, \frac{1}{{(2n+1)}^2\displaystyle {2n \choose n}^2} \cdot \frac{{(n+1)}_k {(n+1)}_k}{{(2n+2)}_k {(2n+2)}_k} \,,\]
which can be verified by straightforward computations. \hfill \qed \\

Next, we define recursively a sequence $p_{\nu}(t)$ $(\nu=1,2,\dots)$ of polynomias in one variable $t$, namely
\begin{eqnarray}
p_1(t) &=& t^2 \,,
\label{5.50} \\
p_2(t) &=& t^4-t^2+t \,,
\label{5.60} \\
p_{\nu}(t) &=& t^2p_{\nu-1}(t) + t(1-t)p_{\nu-2}(t) \qquad (\nu=3,4,\dots).
\label{5.70}
\end{eqnarray}
It is clear that $\deg p_{\nu}=2\nu$, which follows easily by induction for $\nu$ with $\deg p_1=2$
and $\deg p_2=4$. The leading coefficient of $p_{\nu}$ is 1 for $\nu=1,2,\dots$. Let
\[p_{\nu}(t) \,=\, \sum_{\mu=0}^{2\nu} a_{\nu,\mu}t^{\mu} \,.\]
\begin{lemma}
For $\nu \geq 3$ we have
\[p_{\nu}(t) \,=\, t^{2\nu} + (1-\nu) t^{2\nu -2} + \sum_{\mu=2}^{2\nu -3} \big( a_{\nu-1,\mu -2} + 
a_{\nu-2,\mu -1} - a_{\nu-2,\mu-2} \big) t^{\mu} \]
with
\[a_{\nu,\mu} \,=\, a_{\nu-1,\mu -2} + a_{\nu-2,\mu -1} - a_{\nu-2,\mu-2} \qquad (2\leq \mu \leq 2\nu -3)\,.\]
\label{Lem5.1}
\end{lemma}
{\em Proof:\/} \,Using the definition of $p_{\nu}(t)$ from (\ref{5.50}) to (\ref{5.70}) with $\nu \geq 3$,
we obtain
\begin{eqnarray*}
p_{\nu}(t) &=& \sum_{\mu=0}^{2\nu} a_{\nu,\mu}t^{\mu} \,=\,  \sum_{\mu=0}^{2\nu-2} a_{\nu-1,\mu}t^{\mu+2} 
+ \sum_{\mu=0}^{2\nu-4} a_{\nu-2,\mu}t^{\mu+1} - \sum_{\mu=0}^{2\nu-4} a_{\nu-2,\mu}t^{\mu+2}  \\
&=&  \sum_{\mu=2}^{2\nu} a_{\nu-1,\mu-2}t^{\mu} + \sum_{\mu=1}^{2\nu-3} a_{\nu-2,\mu-1}t^{\mu}
- \sum_{\mu=2}^{2\nu-2} a_{\nu-2,\mu-2}t^{\mu}  \\
&=& a_{\nu-1,2\nu-2}t^{2\nu} + a_{\nu-1,2\nu-3}t^{2\nu-1} + a_{\nu-1,2\nu-4}t^{2\nu-2} + a_{\nu-2,0}t -
a_{\nu-2,2\nu-4}t^{2\nu-2} \\
&& + \sum_{\mu=2}^{2\nu-3} \big( a_{\nu-1,\mu-2} + a_{\nu-2,\mu-1} - a_{\nu-2,\mu-2} \big) t^{\mu} \\
&=& t^{2\nu} + (1-\nu)t^{2\nu-2} + \sum_{\mu=2}^{2\nu-3} \big( a_{\nu-1,\mu-2} + a_{\nu-2,\mu-1} - a_{\nu-2,\mu-2} \big)
t^{\mu} \,,
\end{eqnarray*}
since the four identities
\begin{eqnarray*}
a_{\nu-1,2\nu-2} &=& 1\,,\\
a_{\nu-1,2\nu-3} &=& 0 \,,\\
a_{\nu-1,2\nu-4} - a_{\nu-2,2\nu-4} &=& 1-\nu \,,\\
a_{\nu-2,0} &=& 0 
\end{eqnarray*}
follow easily from  (\ref{5.50}) to (\ref{5.70}) by $a_{\nu,2\nu}=1$, $a_{\nu,2\nu-1}=0$, $a_{\nu,2\nu-2}=1-\nu$,
and $a_{\nu,0}=0$ for $\nu\geq 1$. The lemma is proven. \hfill \qed
\begin{thm}
For the sequences $B_2,R_2$ defined above for $\zeta(2)$ we have
\[\E^*\big( B_2,R_2,\zeta(2)\big) \,=\, 1 + \sum_{\nu =1}^{\infty} \,\sum_{\mu =0}^{2\nu} \frac{a_{\nu,\mu}}
{(\nu+1)(\mu+1)} \,.\] 
\label{Thm5.20}
\end{thm}
{\em Proof:\/} \,With (\ref{5.50}) to (\ref{5.70}) we obtain
\begin{eqnarray*} 
&& (1+x^2y^2-xy^2-yx^2) \Big( \,1+\sum_{\nu=1}^{\infty} p_{\nu}(y)x^{ \nu}\,\Big) \\
&=& 1+x^2y^2-xy^2-yx^2 + \sum_{\nu =1}^{\infty} p_{\nu}(y)x^{ \nu} +  \sum_{\nu =1}^{\infty} y^2p_{\nu}
(y)x^{ \nu+2} -  \sum_{\nu =1}^{\infty} y^2p_{\nu}(y)x^{ \nu+1} \\
&& -\,\sum_{\nu =1}^{\infty} yp_{\nu}(y)x^{ \nu+2} \\
&=& 1 + x^2y^2 - xy^2 - yx^2 +  \sum_{\nu =1}^{\infty} p_{\nu}(y)x^{ \nu} - \sum_{\nu =3}^{\infty} 
y(1-y)p_{\nu-2}(y)x^{\nu} - \sum_{\nu=2}^{\infty} y^2p_{\nu-1}(y)x^{\nu} \\
&=& 1 - p_1(y)x - p_2(y)x^2 + y^2p_1(y)x^2  +  \sum_{\nu =1}^{\infty} p_{\nu}(y)x^{ \nu} - 
\sum_{\nu =3}^{\infty} y(1-y)p_{\nu-2}(y)x^{\nu} - \sum_{\nu=2}^{\infty} y^2p_{\nu-1}(y)x^{\nu} \\
&=& 1 +  \sum_{\nu =3}^{\infty} p_{\nu}(y)x^{ \nu} - \sum_{\nu =3}^{\infty} y(1-y)p_{\nu-2}(y)x^{\nu} -
\sum_{\nu=3}^{\infty} y^2p_{\nu-1}(y)x^{\nu} \\
&=& 1 + \sum_{\nu =3}^{\infty} \Big[ \,p_{\nu}(y) - \big( y^2p_{\nu -1}(y)  + y(1-y)p_{\nu -2}(y)\big)\,\Big]
x^{\nu} \\
&=& 1\,.
\end{eqnarray*}
Hence,
\[ \frac{1}{1+x^2y^2-xy^2-yx^2} \,=\, 1 + \sum_{\nu=1}^{\infty} p_{\nu}(y)x^{ \nu} \,=\, 
1 +  \sum_{\nu=1}^{\infty}  \,\Big( \,\sum_{\mu =0}^{2\nu} a_{\nu,\mu}y^{\mu} \,\Big) x^{\nu} .\]
Now the theorem follows from (\ref{5.10}) by two-fold integration with respect to $x$ and $y$. \hfill \qed \\

We can proceed similarly in order to obtain similar results for $\E \big( B_2,R_2,\zeta(2)\big)$ and for 
$\E \big( B_3,R_3,\zeta(3)\big)$. Therefore, we state them without proofs. Again we define recursively a sequence
$q_{\nu}(t)$ $(\nu=1,2,\dots)$ of integer polynomias in one variable
$t$,
\begin{eqnarray*}
q_1(t) &=&2t- t^2 \,,\\
q_2(t) &=& t^4-4t^3+5t^2-t \,,\\
q_{\nu}(t) &=& t(2-t)q_{\nu-1}(t) + t(t-1)q_{\nu-2}(t) \qquad (\nu=3,4,\dots).
\end{eqnarray*}
Let
\[q_{\nu}(t) \,=\, \sum_{\mu=0}^{2\nu} b_{\nu,\mu}t^{\mu} \,.\]
Here, we have
\[\frac{1}{1-x^2y^2-2xy+xy^2+yx^2} \,=\, \sum_{\nu=0}^{\infty} q_{\nu}(y)x^{ \nu} \,=\, 
\sum_{\nu=0}^{\infty}  \,\Big( \,\sum_{\mu =0}^{2\nu} b_{\nu,\mu}y^{\mu} \,\Big) x^{\nu} .\]
\begin{lemma}
For $\nu \geq 3$ we have
\begin{eqnarray*}
q_{\nu}(t) &=&{(-1)}^{\nu} t^{2\nu} +2{(-1)}^{\nu+1}\nu t^{2\nu-1}+{(-1)}^{\nu} (2\nu^2-\nu-1) t^{2\nu -2} \\
&& +\,\sum_{\mu=2}^{2\nu -3} \big( -b_{\nu-1,\mu -2} + 2b_{\nu-1,\mu -1} + b_{\nu-2,\mu-2}  
- b_{\nu-2,\mu-1}\big) t^{\mu}
\end{eqnarray*}
with
\[b_{\nu,\mu} \,=\, -b_{\nu-1,\mu -2} + 2b_{\nu-1,\mu -1} + b_{\nu-2,\mu-2} - b_{\nu-2,\mu-1} \qquad 
(2\leq \mu \leq 2\nu -3)\,.\]
\label{Lem5.2}
\end{lemma}
\begin{thm}
For the sequences $B_2,R_2$ defined above for $\zeta(2)$ we have
\[\E\big( B_2,R_2,\zeta(2)\big) \,=\, 1 + \sum_{\nu =1}^{\infty} \,\sum_{\mu =0}^{2\nu} \frac{b_{\nu,\mu}}
{(\nu+1)(\mu+1)} \,.\] 
\label{Thm5.30}
\end{thm}
The above method can be generalized such that it also works for $\E\big( B_3,R_3,\zeta(3)\big)$. Let
\begin{eqnarray*}
r_0(x,y) &=& 1 \,,\nonumber \\
r_1(x,y) &=& x^2y^2-xy^2-x^2y+1 \,,\nonumber \\
r_2(x,y) &=& x^4y^4-2x^3y^4-2x^4y^3+x^2y^4+x^4y^2+2x^3y^3+x^2y^2-xy^2-x^2y-xy+1 \,,
\nonumber \\
r_{\nu}(x,y) &=& \big( x^2y^2-xy^2-x^2y+1 \big) r_{\nu -1}(x,y) - \big( x^2y^2-xy^2-x^2y+xy \big) 
r_{\nu -2}(x,y) \,,
\end{eqnarray*}
where $\nu \geq 3$. Setting
\[r_{\nu}(x,y) \,=\, \sum_{\mu_1=0}^{2\nu}\,\sum_{\mu_2=0}^{2\nu}  c_{\nu,\mu_1,\mu_2} x^{\mu_1}
y^{\mu_2} \,,\]
it turns out that
\[\frac{1}{1+x^2y^2w^2-xy^2w^2-x^2yw^2-x^2y^2w+xyw^2+xy^2w+x^2yw-w} 
\,=\, \sum_{\nu =0}^{\infty} r_{\nu}(x,y)\cdot w^{\nu} \,.\]
Then, (\ref{5.30}) underlies the following result.
\begin{thm}
For the sequences $B_3,R_3$ defined above for $\zeta(3)$ we have 
\[\E\big( B_3,R_3,\zeta(3)\big) \,=\, \frac{1}{2} + \frac{1}{2} \sum_{\nu =1}^{\infty} \, \sum_{\mu_1=0}^{2\nu}\,
\sum_{\mu_2=0}^{2\nu} \frac{c_{\nu,\mu_1,\mu_2}}{(\nu +1)(\mu_1 +1)(\mu_2 +1)} \,.\]
\label{Thm5.40}
\end{thm}
{\bf 4.)} As mentionned at the beginning of Section~\ref{Sec5}, some connections between rational coefficients involved in computing the error sums for Ap\'ery's continued fraction and well-known integer sequences may be noticed. Furthermore some unproved identities have been empirically found for such coefficients.

The Theorems~\ref{5.20} and~\ref{5.30} rely on two triangles of integer coefficients, namely~$a_{\nu,\mu}$ for the former and~$b_{\nu,\mu}$ for the latter. Both can be expressed by binomial sums as follows.
\[
\begin{array}{lcl}
a_{\nu,\mu}&=&\displaystyle
  \sum_{k=0}^\nu\sum_{i=0}^k
  \left(-1\right)^{\nu+k}
  {{\mu-k}\choose{2\mu-\nu-k-i}}
  {{\mu-k}\choose i}
  {{\mu-i}\choose{k-i}} \,, \\[20pt]
b_{\nu,\mu}&=&\displaystyle
  \sum_{k=0}^\nu\sum_{i=0}^k
  \left(-1\right)^{\nu+\mu}
  {{\mu-k}\choose{2\mu-\nu-k-i}}
  {{\mu-k}\choose i}
  {{\mu-i}\choose{k-i}} \,,
\end{array}
\]
which both lead to non-recurrent formulas for the error sums as quadruple sums.

Several basic properties concerning the coefficients $a_{\nu,\mu}$ and~$b_{\nu,\mu}$ can be noticed, including
\[
\displaystyle\sum_{\mu=0}^{2\nu}a_{\mu,\nu}=1
\qquad\textrm{and}\qquad
\displaystyle\sum_{\mu=0}^{2\nu}b_{\mu,\nu}=1
\qquad\qquad\left(\nu\in\mathbb{N}\right)
\]
and
\[
a_{\nu,\mu}=a_{\mu,\nu}\qquad\textrm{and}\qquad b_{\nu,\mu}=b_{\mu,\nu} \,.
\]

More unproved identities come from the theory of generating functions. Both coefficients~$a_{\nu,\mu}$ and~$b_{\nu,\mu}$ seem to be the coefficients of degree~$2\nu-\mu$ in the MacLaurin series expansion of
\[
\left\{
\begin{array}{l@{\qquad}l}
\displaystyle\frac{
  \left(
    \displaystyle\frac
      { x^2 + 1 - \sqrt{x^4-4x^3+2x^2+1} }
      { 2x^3}
  \right)^{\mu-\nu}
}{
  \sqrt{x^4-4x^3+2x^2+1}
}&\textrm{for $a_{\nu,\mu}$}
\\[20pt]
\displaystyle\frac{
  \left(
    \displaystyle\frac
      { x^2 + 2x - 1 + \sqrt{x^4+2x^2-4x+1} }
      { 2x^3}
  \right)^{\mu-\nu}
}{
  \sqrt{x^4+2x^2-4x+1}
}&\textrm{for $b_{\nu,\mu}$}
\end{array}\right. \quad .
\]
These generating functions actually allow to build the triangles of coefficients~$a_{\nu,\mu}$ and~$b_{\nu,\mu}$ by diagonals rather than by rows.

Summing these coefficients by rows according to Theorems~\ref{5.20} and~\ref{5.30}, the results can be easely achieved by
applying the following unproved recursive identities.
\[
\left\{
\begin{array}{l}
\alpha_0=1,\qquad\alpha_1=1/3,\qquad\alpha_2=11/30,\qquad\alpha_3=17/70\\[12pt]
\begin{array}{lclcl}
\alpha_{n}& =& \displaystyle\frac{ 4n-1} { 2n+1 } \alpha_{n-1}
&-& \displaystyle\frac{2n-2 }{ 2n+1} \alpha_{n-2}\\[12pt]
&-&
\displaystyle\frac{n-1 } {4n+2} \alpha_{n-3}
&+& \displaystyle\frac{ n-2 }{4n+2} \alpha_{n-4}
\end{array}
\end{array}\right. \qquad 
\]
and
\[
\left\{
\begin{array}{l}
\beta_0=1,\qquad\beta_1=2/3,\qquad\beta_2=11/30,\qquad\beta_3=47/210\\[12pt]
\begin{array}{lclcl}
\beta_{n}&=& \displaystyle\frac{ 6n-1} { 2n+1 } \beta_{n-1}
&-& \displaystyle\frac{6n-5 }{ 2n+1} \beta_{n-2}\\[12pt]
&+&
 \displaystyle\frac{5n-7 } {4n+2} \beta_{n-3}
&-& \displaystyle\frac{n-2 }{4n+2} \beta_{n-4}
\end{array}
\end{array}\right. \quad ,
\]
where
\[
\alpha_\nu=\displaystyle\sum_{\mu=0}^{2\nu}\displaystyle\frac{a_{\nu,\mu}}{\mu+1}
\qquad\textrm{and}\qquad
\beta_\nu=\displaystyle\sum_{\mu=0}^{2\nu}\displaystyle\frac{b_{\nu,\mu}}{\mu+1} \,.
\]

The special case~$b_{n,n}$, which may be called the main diagonal in the triangle of coefficients~$b_{\nu,\mu}$, leads to the following simplifications. We have
\[
b_{n,n}=
\displaystyle\sum_{k=0}^n\sum_{i=0}^k
  \displaystyle
  {{n-k}\choose i}^2
  {{n-i}\choose{k-i}} \,, 
\]
where the generating function of the $b_{n,n}$ is given by $1/\sqrt{x^4+2x^2-4x+1}$.
This is the sequence A108626 from the {\em On-Line Encyclopedia of Integer Sequences}. This sequence gives the antidiagonal sums of the square array~A108625 itself known to be highly related to the constant~$\zeta(2)$.

$b_{\nu,\mu}$ is defined recursively by
\[
b_{\nu,\mu}=2b_{\nu-1,\mu-1}-b_{\nu-1,\mu-2}+b_{\nu-2,\mu-2}-b_{\nu-2,\mu-1} \,.
\]
Assuming $b_{n,n+1}=b_{n+1,n}$ (unproved), a new recursive identity can be given concerning A108626:
\[
\begin{split}
\texttt{A108626}\left(n+2\right)-2\times\texttt{A108626}\left(n+1\right)-\texttt{A108626}\left(n\right)\\=2\displaystyle\sum_{k=0}^n\sum_{i=0}^k{{n-k+1}\choose{i-1}}{{n-k+1}\choose i}{{n-i+1}\choose{k-i}} \,.
\end{split}
\]
The previous relation actually happens to be the simplest case from a more general sequence of recurrence relations of order~$2d$ given by:
\[
\displaystyle\sum_{k=0}^{2d} c_{k}\texttt{A108626}\left(n+k\right)
  = \left(-1\right)^d
  \displaystyle\sum_{k=0}^n\sum_{i=0}^k
    {{n-k}\choose{d+i}}{{n-k}\choose i}{{n-i}\choose{k-i}} \,,
\]
where the numbers $c_k$ are coefficients of order~$2d-k$ in the characteristic polynomial
\[
\displaystyle\frac{1}{2^d}\sum_{i=0}^{\left\lfloor\frac{d}{2}\right\rfloor}
  {d\choose{2i}}\left(x^4+2x^2-4x+1\right)^i\left(x^2+2x-1\right)^{d-2i} \,.
\]

These recurrence relations, as well as similar ones related to the coefficients $a_{\nu,\mu}$, can be written as new generating functions, the diagonal of order~$d$ being made from the coefficients of terms with positive powers in
\[
\left\{
\begin{array}{l@{\qquad}l}
\displaystyle\frac{
  \displaystyle\sum_{k=0}^{\left\lfloor\frac{d}{2}\right\rfloor}
    {d\choose{2k}}\left(x^4-4x^3+2x^2+1\right)^k\left(x^2+1\right)^{d-2k}
}{ \left(2x^3\right)^d \sqrt{x^4-4x^3+2x^2+1} }
&\textrm{for $a_{n,n+d}$}
\\[20pt]
\displaystyle\frac{
  \displaystyle\sum_{k=0}^{\left\lfloor\frac{d}{2}\right\rfloor}
    {d\choose{2k}}\left(x^4+2x^2-4x+1\right)^k\left(x^2+2x-1\right)^{d-2k}
}{ \left(2x^3\right)^d \sqrt{x^4+2x^2-4x+1} }
&\textrm{for $b_{n,n+d}$}
\end{array}\right. \quad .
\]

%


\[\]
List of OEIS sequence numbers. \\

A001850, \\
A003417, \\
A005258, \\
A005259, \\
A051451, \\
A108626, \\
A108626.

\end{document}